\documentclass{conm-p-l}

\theoremstyle{definition}

\theoremstyle{remark}

\numberwithin{equation}{section}

\usepackage{amssymb}
\newcommand\CC{\mathbb{C}}
\newcommand\RR{\mathbb{R}}
\newcommand\ZZ{\mathbb{Z}}
\newcommand\FSD{\mathcal{D}}
\newcommand\al\alpha
\newcommand\be\beta
\newcommand\ga\gamma
\newcommand\tha\theta
\newcommand\la\lambda
\newcommand\Om\Omega
\newcommand\half{\frac12}
\newcommand\iy\infty
\newcommand\pa\partial
\newcommand\thalf{\textstyle\half}

\newcommand\RHS{right-hand side}
\newcommand\const{\mbox{const.}\,}
\newcommand\Znonneg{\ZZ_{\ge0}}
\newcommand{\hyp}[5]{\,\mbox{}_{#1}F_{#2}\left(
  \genfrac{}{}{0pt}{}{#3}{#4};#5\right)}
\begin{document}
\title[Lowering and raising operators]
{Lowering and raising operators for some special\\
orthogonal polynomials}
%
\author{Tom H. Koornwinder}
\address{Korteweg-de Vries Institute, University of Amsterdam,
Plantage Muidergracht 24, 1018 TV Amsterdam, The Netherlands}
\email{thk@science.uva.nl}
%
\subjclass{Primary 33C45, 33D45, 33C52}
\date{February 21, 2005.}
\dedicatory{This paper is dedicated to Ian Macdonald on the occasion of his
75th birthday.}
\keywords{lowering and raising operators, continuous q-ultraspherical
polynomials, ultraspherical polynomials, Jacobi polynomials,
Jacobi polynomials for root system $BC_2$}
\begin{abstract}
This paper discusses operators lowering or raising the degree
but preserving the parameters of special orthogonal
polynomials. Results for one-variable classical (q-)orthogonal
polynomials are surveyed. For Jacobi polynomials associated with
root system $BC_2$ a new pair of lowering and raising operators
is obtained.
\end{abstract}
\maketitle
%
\section{Introduction} 
\label{sec:05}
Kirillov and Noumi \cite{10} gave explicit $q$-difference lowering and raising
operators for $A_{n-1}$ type Macdonald polynomials
$J_\la(x;q,t)=c_\la(q,t) P_\la(x;q,t)$ (see \cite[(VI.8.3)]{11}).
These operators don't change the parameter $t$, they only lower or
raise $\la$. This is quite different from
Opdam's \cite{05} shift operators acting on Jacobi polynomials associated
with root systems, which do change parameters.
In \cite[Remark 5.3]{10} the interesting problem is
mentioned to find lowering and raising operators for
Macdonald-Koornwinder polynomials (see \cite{12}). As far as we know,
such operators have not yet been given in literature until now,
and neither in the corresponding $q=1$ case of $BC_n$ type
Jacobi polynomials (see \cite{03}, \cite{04}).

The present paper makes only a minor step in the direction of a general
answer to the problem raised in \cite[Remark 5.3]{10}. In the $BC_2$,
$q=1$ case, for parameters $(\al,\be,\ga)$ with $\al=\be$,
and for partition $\la=(n,0)$ of length 1, a raising and lowering operator
in explicit form are obtained (sections \ref{sec:03} and
\ref{sec:04}).
In the earlier sections \ref{sec:06}, \ref{sec:02}, \ref{sec:01}
lowering and raising
formulas in the rank 1 case
(continuous $q$-ultraspherical, ultraspherical and Jacobi polynomials,
but not yet Askey-Wilson polynomials) are discussed.
Some formulas scattered in the literature are brought here together,
and also the explicit specialization of the Kirillov-Noumi operators
to the $A_1$ case is given. Some of the formulas in
sections \ref{sec:06}--\ref{sec:01}
may be new in one-variable special function theory.
\subsubsection*{Acknowledgement}
I thank Hechun Zhang (Tsinghua University, Bejing) for helpful
discussions in an early stage of writing this manuscript.
\section{Lowering and raising continuous $q$-ultraspherical polynomials}
\label{sec:06}
See \cite{14} or \cite[\S3.10.1]{13} for the definition of the
\textit{continuous $q$-ultraspherical polynomials} $C_n(x;t|q)$.
A pleasant explicit formula for them is as a finite Fourier series:
\begin{equation}
C_n(\cos\tha;t|q)=\sum_{k=0}^n \frac{(t;q)_k\,(t;q)_{n-k}}
{(q;q)_k\,(q;q)_{n-k}}\,e^{i(n-2k)\tha}.
\label{eq:06.15}
\end{equation}
(Throughout see \cite{15} for definition of ($q$-)hypergeometric series
and ($q$-)Pochhammer symbols.)
The $A_1$ type Macdonald polynomials can be expressed in terms
of continuous $q$-ultraspherical polynomials:
\begin{equation}
P_{m,n}(x,y;q,t)=
\frac{(q;q)_{m-n}}{(t;q)_{m-n}}\,(xy)^{\half(m+n)}\,
C_{m-n}\biggl(\frac{x+y}{2(xy)^\half};t|q\biggr)\quad
(m\ge n\ge0).
\end{equation}
For the renormalized polynomials we have
\begin{multline*}
J_{m,n}(x,y;q,t)=c_{m,n}(q,t) P_{m,n}(x,y;q,t),\\
\mbox{where}\quad c_{m,n}(q,t)=(t^2 q^{m-n};q)_n\,(t;q)_{m-n}\,(t;q)_n,
\quad\mbox{so:}
\end{multline*}
\begin{multline}
J_{m,n}(x,y;q,t)=(t^2 q^{m-n};q)_n\,(t;q)_n\,(q;q)_{m-n}\\
\times(xy)^{\half(m+n)}\,
C_{m-n}\biggl(\frac{x+y}{2(xy)^\half};t|q\biggr)
\quad(m\ge n\ge0).
\end{multline}
In particular, if $(m,n)$ is replaced by $(n,0)$, and
if we use \eqref{eq:06.15}:
\begin{align}
J_{n,0}(x,y;q,t)&=(q;q)_n\,(xy)^{\half n}\,
C_n\biggl(\frac{x+y}{2(xy)^\half};t|q\biggr)
\label{eq:06.03}\\
&=(q;q)_n\,\sum_{k=0}^n
\frac{(t;q)_k\,(t;q)_{n-k}}{(q;q)_k\,(q;q)_{n-k}}\,x^{n-k}y^k.
\label{eq:06.14}
\end{align}

Let us consider the Kirillov-Noumi operators \cite{10} in the two-variable
case while acting on the polynomials $J_{n,0}$, raising or lowering $n$.
There are two variants $K_1^+$ and $K_1^-$ of the raising operator,
and two variants $M_1^+$ and $M_1^-$ of the lowering operator.
The expressions for these operators will involve the operators $T_{q,x}$ and
$T_{q,y}$:
\begin{equation*}
(T_{q,x}f)(x,y):=f(qx,y),\quad
(T_{q,y}f)(x,y):=f(x,qy).
\end{equation*}
Then $K_1^\pm$ and $M_1^\pm$ and their actions are given by:
\begin{align*}
K_1^+&:=x\left(1-\frac{tx-y}{x-y}\,T_{q,x}\right)
+y\left(1-\frac{ty-x}{y-x}\,T_{q,y}\right),\\
K_1^-&:=x\left(-t T_{q,x} T_{q,y}+\frac{x-ty}{x-y}\,T_{q,y}\right)+
y\left(-t T_{q,x} T_{q,y}+\frac{y-tx}{y-x}\,T_{q,x}\right),
\end{align*}
\begin{equation}
K_1^\pm\,J_{n,0}(x,y;q,t)=J_{n+1,0}(x,y;q,t);
\label{eq:06.01}
\end{equation}
\begin{align*}
M_1^+&:=\frac1x\left(1-\frac{tx-y}{x-y}\,T_{q,x}\right)
+\frac1y\left(1-\frac{ty-x}{y-x}\,T_{q,y}\right),\\
M_1^-&:=\frac1x\left(-t T_{q,x} T_{q,y}+\frac{x-ty}{x-y}\,T_{q,y}\right)+
\frac1y\left(-t T_{q,x} T_{q,y}+\frac{y-tx}{y-x}\,T_{q,x}\right),
\end{align*}
\begin{equation}
M_1^\pm\,J_{n,0}(x,y;q,t)=(1-q^n)(1-t^2q^{n-1})\,J_{n-1,0}(x,y;q,t).
\label{eq:06.02}
\end{equation}
For the moment I will take \eqref{eq:06.01} and \eqref{eq:06.02} for
granted from \cite{10}. At the end of this section I will prove these
formulas independently.

We can write \eqref{eq:06.01} and \eqref{eq:06.02} more explicitly
by substitution of
\eqref{eq:06.03} and the explicit expressions for the operators
$K_1^\pm$ and $M_1^\pm$. Then put $x=z$, $y=z^{-1}$. In the resulting
formulas it is
convenient to assume $t,q$ fixed and to use the notation
\begin{equation*}
C_n[z]:=C_n(\thalf(z+z^{-1});t|q),\qquad
TC_n[z]:=C_n[q^\half z],\qquad
T^{-1}C_n[z]:=C_n[q^{-\half} z].
\end{equation*}
We obtain:
\begin{multline}
-\frac{tz^2-1}{z-z^{-1}}\,TC_n[z]+
\frac{tz^{-2}-1}{z-z^{-1}}\,T^{-1}C_n[z]+
q^{-\half n}(z+z^{-1})\,C_n[z]\\
=(q^{-\half n}-q^{\half n+1})\,C_{n+1}[z],
\label{eq:06.04}
\end{multline}
\begin{multline}
-\frac{z^{-2}-t}{z-z^{-1}}\,TC_n[z]+
\frac{z^2-t}{z-z^{-1}}\,T^{-1}C_n[z]-
q^{\half n}t(z+z^{-1})\,C_n[z]\\
=(q^{-\half n}-q^{\half n+1})\,C_{n+1}[z],
\label{eq:06.05}
\end{multline}
\begin{multline}
-\frac{t-z^{-2}}{z-z^{-1}}\,TC_n[z]+
\frac{t-z^2}{z-z^{-1}}\,T^{-1}C_n[z]+
q^{-\half n}(z+z^{-1})\,C_n[z]\\
=(q^{-\half n}-t^2q^{\half n-1})\,C_{n-1}[z],
\label{eq:06.06}
\end{multline}
\begin{multline}
-\frac{1-tz^2}{z-z^{-1}}\,TC_n[z]+
\frac{1-tz^{-2}}{z-z^{-1}}\,T^{-1}C_n[z]-
q^{\half n}t(z+z^{-1})\,C_n[z]\\
=(q^{-\half n}-t^2q^{\half n-1})\,C_{n-1}[z].
\label{eq:06.07}
\end{multline}

If we subtract \eqref{eq:06.05} from \eqref{eq:06.04} or
\eqref{eq:06.07} from \eqref{eq:06.06}, and if we divide the resulting
second order $q^\half$-difference formula by a suitable factor which
all terms have in common, then we obtain
\begin{equation}
\frac{1-tz^2}{1-z^2}\,TC_n[z]+
\frac{1-tz^{-2}}{1-z^{-2}}\,T^{-1}C_n[z]-
(q^{-\half n}+q^{\half n} t)\,C_n[z]=0,
\label{eq:06.08}
\end{equation}
which is a known formula for continuous
$q$-ultraspherical polynomials.
Indeed, rewrite continuous $q$-ultraspherical polynomials
as Askey-Wilson polynomials by
\begin{multline*}
C_n(x;q^{\al+\half}|q)=
\const P_n^{(\al,\al)}(x;q^\half)\\
=\const p_n(x;q^{\frac14},q^{\half\al+\frac14},-q^{\half\al+\frac14},
-q^{\frac14}|q^\half),
\end{multline*}
(see \cite[(7.5.34), (7.5.25), (7.5.1)]{15}), and then use the
second order $q$-difference formula
\cite[(3.1.7)]{13}.
Thus \eqref{eq:06.04} is equivalent with \eqref{eq:06.05} modulo
\eqref{eq:06.08}, and similarly for \eqref{eq:06.06} and \eqref{eq:06.07}.

In addition to the operators $K^\pm$ and $L^\pm$ we introduce the operators
$A$ and $\Om$ given by:
\begin{equation*}
A:=T_{q,x}\,T_{q,y},\qquad
\Om:=\frac1{xy}\,
\left(1-\frac{tx-y}{x-y}\,T_{q,x}
-\frac{x-ty}{x-y}T_{q,y}+t\,T_{q,x}\,T_{q,y}
\right).
\end{equation*}
Since
\begin{equation*}
K_1^+-K_1^-=xy(x+y)\,\Om,\qquad
M_1^+-M_1^-=(x+y)\,\Om,
\end{equation*}
we will no longer consider $K_1^-$ and $M_1^-$, but we will concentrate
on $K_1^+$, $M_1^+$, $\Om$ and $A$. We can derive relations
\begin{equation}
A\,\Om=q^{-2}\,\Om\,A,
\qquad
A\,K_1^+=q\,K_1^+\,A,
\qquad
A\,M_1^+=\,q^{-1}\,M_1^+\,A,
\label{eq:06:16}
\end{equation}
\begin{equation*}
\Om\,K_1^+=q^2\,K_1^+\,\Om+(1-q)^2(x+y)\,\Om,
\qquad
\Om\,M_1^+=M_1^+\,\Om,
\end{equation*}
\begin{equation*}
q^2\,K_1^+\,M_1^+-M_1^+\,K_1^+=(q^2-1)+(1-q)(q+t^2)\,A
+(q^2-1)(x^2+xy+y^2)\,\Om,
\end{equation*}

A straightforward computation shows that the operators
$(x-y)\,\Om$, $(x-y)\,K_1^+$ and $(x-y)\,M_1^+$ send $x^my^n+x^ny^m$
to an antisymmetric polynomial in $x$ and $y$.
Therefore $\Om$, $K_1^\pm$, $M_1^\pm$ and (clearly) $A$ act on the space
of symmetric polynomials in $x$ and $y$. By the above relations,
the operators $K_1^+$, $M_1^+$ and $A$ also act on the subspace
of symmetric polynomials annihilated by $\Om$. The operators $K_1^+$,
$M_1^+$ and $A$ restricted to this subspace satisfy the relations
\begin{equation}
q^2\,K_1^+\,M_1^+-M_1^+\,K_1^+=(q^2-1)+(1-q)(q+t^2)\,A,
\label{eq:06.09}
\end{equation}
\begin{equation}
A\,K_1^+=q\,K_1^+\,A,
\qquad
A\,M^+=\,q^{-1}\,M_1^+\,A.
\label{eq:06.10}
\end{equation}

It does not seem that the relations \eqref{eq:06.09} and
\eqref{eq:06.10} are equivalent to the familiar relations
\begin{equation}
EF-FE=\frac{K-K^{-1}}{q-q^{-1}}\,,\qquad
KE=q^2 EK,\qquad
KF=q^{-2}FK,
\label{eq:06.11}
\end{equation}
for the
generators of $U_q(sl(2))$ as given, for instance, in
\cite[Definition VI.1.1]{16}.
Indeed, relations \eqref{eq:06.09} and \eqref{eq:06.10}
take after rescaling of the generators the form
\begin{equation}
q^2\,K_1^+\,M_1^+-M_1^+\,K_1^+=A-1,
\qquad
A\,K_1^+=q\,K_1^+\,A,
\qquad
A\,M^+=\,q^{-1}\,M_1^+\,A.
\label{eq:06.12}
\end{equation}
while relations \eqref{eq:06.11}, after substitution of $\tilde E:=EK$
and after rescaling,
become
\begin{equation}
q^2\tilde EF-F\tilde E=K^2-1,\qquad
K\tilde E=q^2 \tilde EK,\qquad
KF=q^{-2}FK.
\label{eq:06.13}
\end{equation}
Relations \eqref{eq:06.12} would match with relations \eqref{eq:06.13}
if the first relation in \eqref{eq:06.12} would have been
$q^\half K_1^+\,M_1^+-M_1^+\,K_1^+=A-1$.

Now I will give the promised independent proof of
\eqref{eq:06.01} and \eqref{eq:06.02}.
Because of the first relation in \eqref{eq:06:16},
the symmetric polynomials annihilated by $\Om$ have a basis of
homogeneous polynomials.
From $\Om\left(\sum_{k=0}^n c_k(x^{n-k}y^k+x^ky^{n-k})\right)=0$ with
$c_k=c_{n-k}$
one derives a recurrence relation for the $c_k$ which, on comparison
with \eqref{eq:06.14}, shows that the polynomials
$J_{n,0}(x,y;q,t)$ ($n\in\Znonneg$) given by \eqref{eq:06.03} span the
space of symmetric polynomials annihilated by $\Om$.
Thus, because $K_1^+$ resp.\ $M_1^+$ raise resp.\ lower the degree of
a homogeneous symmetric polynomial by 1, we find
\eqref{eq:06.01} for $K_1^+$ and \eqref{eq:06.02} for $M_1^+$
up to a constant factor. These constant factors are then obtained
by comparing terms of highest degree on the left and on the right.

\section{Lowering and raising ultraspherical polynomials}
\label{sec:02}
%
The \emph{ultraspherical polynomials} $C_n^{(\la)}(x)$
(see for instance \cite{01} or \cite{13})
can be obtained from
\eqref{eq:06.15} by putting $t=q^\la$ and letting $q\uparrow1$:
\begin{equation*}
C_n^{(\la)}(\cos\tha)=\sum_{k=0}^n \frac{(\la)_k\,(\la)_{n-k}}
{k!\,(n-k)!}\,e^{i(n-2k)\tha}.
\end{equation*}
They are special cases of Jacobi polynomials (see \eqref{eq:01.13}):
\begin{equation*}
C_n^{(\la)}(x)=\frac{(2\la)_n}{(\la+\thalf)_n}\,P_n^{(\la-\half,\la-\half)}(x).
\end{equation*}

In principle, one could take limits for $q\uparrow1$
of all formulas in \S\ref{sec:06}, but I prefer to present the
results on lowering and raising operators for ultraspherical polynomials,
more classical than the results in the $q$-case, here independently from
\S\ref{sec:06}. Lowering and raising formulas can be given in three
different forms.

The first form (see \cite[10.9(15)]{01}) is:
\begin{align}
\left((1-x^2)\frac d{dx}+nx\right) C_n^{(\la)}(x)=&
(n+2\la-1)\,C_{n-1}^{(\la)}(x),
\label{eq:02.06}\\
\left((1-x^2)\frac d{dx}-(n-1+2\la)x\right) C_{n-1}^{(\la)}(x)=&
-n\,C_n^{(\la)}(x).
\label{eq:02.07}
\end{align}
Note that substitution of \eqref{eq:02.06} into \eqref{eq:02.07}
causes $n$ to drop out from the terms with derivatives. There results
$(1-x^2)$ times the second order differential equation
\cite[10.9(14)]{01} for $C_n^{(\la)}(x)$.
Also, if $n$ is replaced by $n+1$ in \eqref{eq:02.07} and if
the term with first derivative is eliminated from the
resulting equation together with \eqref{eq:02.06}, then we obtain
the three-term recurrence relation \cite[10.9(13)]{01} for
$C_n^{(\la)}(x)$.

We get a second form of lowering and raising formulas by rewriting
\eqref{eq:02.06} and \eqref{eq:02.07} into an equivalent form:
\begin{equation}
\frac d{dx}\left((1+x^2)^{\half n}\,
C_n^{(\la)}\left(\frac x{\sqrt{1+x^2}}\right) \right)
=(n+2\la-1)\,(1+x^2)^{\half(n-1)}\,
C_{n-1}^{(\la)}\left(\frac x{\sqrt{1+x^2}}\right),
\label{eq:02.01}
\end{equation}
\begin{equation}
\frac d{dx}\left((1+x^2)^{-\half(n-1)-\la}\,
C_{n-1}^{(\la)}\left(\frac x{\sqrt{1+x^2}}\right) \right)
=-n\,(1+x^2)^{-\half n-\la}\,
C_n^{(\la)}\left(\frac x{\sqrt{1+x^2}}\right).
\label{eq:02.02}
\end{equation}
Iteration of \eqref{eq:02.02} yields the Rodrigues type formula
\begin{equation}
C_n^{(\la)}\left(\frac x{\sqrt{1+x^2}}\right)=
\frac{(-1)^n}{n!}\,(1+x^2)^{\half n+\la}\,
\frac{d^n}{dx^n}\left((1+x^2)^{-\la}\right).
\label{eq:02.03}
\end{equation}
A formula equivalent to \eqref{eq:02.03} (by analytic continuation)
is given in \cite[10.9(37)]{01}, where the formula is ascribed
to F. Tricomi, Ann.\ Mat.\ Pura Appl.\ (4) \textbf{28} (1949), 283--300
(but I could not find the formula there).

Tranformation of the generating function \cite[10.9(29)]{01} for
$C_n^{(\la)}(x)$ yields
\begin{multline*}
\frac 1{(1+(x-z)^2)^\la}=
\sum_{n=0}^\iy (1+x^2)^{-\half n-\la}\,
C_n^{(\la)}\left(\frac x{\sqrt{1+x^2}}\right)z^n\\
(z\in\CC,\;x\in\RR,\;|z|<\sqrt{1+x^2}).
\end{multline*}
Then \eqref{eq:02.03} follows by considering Taylor coefficients
in the above formula.

We obtain a third form of lowering and raising operators by
rewriting \eqref{eq:02.06} and \eqref{eq:02.07}
in an equivalent form as follows:
\begin{multline}
\left(\frac\pa{\pa x}+\frac\pa{\pa y}\right)
\left((xy)^{\half n}\,
C_n^{(\la)}\left(\thalf\bigl((x/y)^\half+(y/x)^\half\bigr)\right)\right)\\
=(n+2\la-1)\,(xy)^{(\half(n-1)}\,
C_{n-1}^{(\la)}\left(\thalf\bigl((x/y)^\half+(y/x)^\half\bigr)\right),
\label{eq:02.10}
\end{multline}
\begin{multline}
\left(x^2\frac\pa{\pa x}+y^2\frac\pa{\pa y}+\la(x+y)\right)
\left((xy)^{\half(n-1)}\,
C_{n-1}^{(\la)}\left(\thalf\bigl((x/y)^\half+(y/x)^\half\bigr)\right)\right)\\
=n\,(xy)^{\half n}\,
C_n^{(\la)}\left(\thalf\bigl((x/y)^\half+(y/x)^\half\bigr)\right).
\label{eq:02.11}
\end{multline}
Iteration of \eqref{eq:02.11} yields the Rodrigues type formula
\begin{equation}
(xy)^{\half n}\,
C_n^{(\la)}\left(\thalf\bigl((x/y)^\half+(y/x)^\half\bigr)\right)
=\frac 1{n!}\left(x^2\frac\pa{\pa x}+y^2\frac\pa{\pa y}+\la(x+y)\right)^n(1).
\label{eq:02.16}
\end{equation}

Formulas \eqref{eq:02.10} and \eqref{eq:02.11}
may be rewritten in terms of the following
special \emph{Jack polynomials} in
two variables:
\begin{equation*}
J_{n,0}^{1/\la}(x,y)=\frac{(\la)_n}{\la^n}P_{n,0}^{1/\la}(x,y)=
\frac{n!}{\la^n}\,(xy)^{\half n}\,
C_n^\la\left(\thalf\bigl((x/y)^\half+(y/x)^\half\bigr)\right)
\end{equation*}
(see also \eqref{eq:03.06}, \eqref{eq:03.07}).
Thus \eqref{eq:02.10} and \eqref{eq:02.11} can be seen to be special cases
of formualas (5.14) resp.\ (2.16) in \cite{10} (there put
$n=2$, $m=1$).

Formulas \eqref{eq:02.10} and \eqref{eq:02.11} are realizations of a
representation of the Lie algebra $sl(2)$. Indeed, put
\begin{equation*}
H:=2\left( x\frac\pa{\pa x}+y\frac\pa{\pa y}+\la\right),\;\;
E:=x^2\frac\pa{\pa x}+y^2\frac\pa{\pa y}+\la(x+y),\;\;
F:=-\left(\frac\pa{\pa x}+\frac\pa{\pa y}\right),
\end{equation*}
\begin{equation*}
f_n(x,y):=(xy)^{\half n}\,
C_n^{(\la)}\left(\thalf\bigl((x/y)^\half+(y/x)^\half\bigr)\right).
\end{equation*}
Then
\begin{align*}
&[H,E]=2E,\quad
[H,F]=-2F,\quad
[E,F]=H,\\
&Ff_n=-(n+2\la-1) f_{n-1},\quad
Ef_{n-1}=n f_n,\quad
Hf_n=2(n+\la) f_n.
\end{align*}
\section{Lowering and raising Jacobi polynomials}
\label{sec:01}
%
%
\emph{Jacobi polynomials} $P_n^{(\al,\be)}(x)$ and their normalized version
$R_n^{(\al,\be)}(x)$ (see for instance \cite{01} or \cite{13})
can be defined in terms
of hypergeometric functions by
\begin{equation}
\begin{split}
P_n^{(\al,\be)}(x)&=
\frac{(\al+1)_n}{n!}\,R_n^{(\al,\be)}(x)\\
&=\frac{(\al+1)_n}{n!}\hyp21{-n,n+\al+\be+1}{\al+1}{\thalf(1-x)}.
\end{split}
\label{eq:01.13}
\end{equation}
They satisfy the symmetry
\begin{equation}
P_n^{(\al,\be)}(-x)=(-1)^n P_n^{(\be,\al)}(x),
\label{eq:01.11}
\end{equation}
and the second order differential equation
(see \cite[10.8(14)]{01})
\begin{equation}
\left(
(1-x^2)\frac{d^2}{dx^2}+\bigl(\be-\al-(\al+\be+2)x\bigr)\frac d{dx}+
n(n+\al+\be+1)
\right)
P_n^{(\al,\be)}(x)=0.
\label{eq:01.04}
\end{equation}

As in \S\ref{sec:02}, lowering and raising formulas can be given in
three different forms. I~start with analogues of
\eqref{eq:02.01} and \eqref{eq:02.02} (the second form in \S\ref{sec:02}):
\begin{multline}
\left(\frac{d^2}{dx^2}+\frac{2\al+1}x\,\frac d{dx}\right)
\left((1+x^2)^n\,P_n^{(\al,\be)}\left(\frac{1-x^2}{1+x^2}\right)\right)\\
=-4(n+\al)(n+\be)\,(1+x^2)^{n-1}
P_{n-1}^{(\al,\be)}\left(\frac{1-x^2}{1+x^2}\right),
\label{eq:01.01}
\end{multline}
\begin{multline}
\left(\frac{d^2}{dx^2}+\frac{2\al+1}x\,\frac d{dx}\right)
\left((1+x^2)^{-n-\al-\be}\,P_{n-1}^{(\al,\be)}\left(\frac{1-x^2}{1+x^2}
\right)\right)\\
=-4n(n+\al+\be)\,(1+x^2)^{-n-\al-\be-1}
P_n^{(\al,\be)}\left(\frac{1-x^2}{1+x^2}\right).
\label{eq:01.06}
\end{multline}
These formulas were first given in \cite[(2.10), (2.11)]{02}.

From \eqref{eq:01.01} and \eqref{eq:01.06} one can derive analogues
of \eqref{eq:02.06} and \eqref{eq:02.07} (the first form in \S\ref{sec:02}):
\begin{multline}
\left(
(2n+\al+\be)(1-x^2)\frac d{dx}+n\bigl((2n+\al+\be)x+\be-\al\bigr)
\right)
P_n^{(\al,\be)}(x)\\
=2(n+\al)(n+\be)\,P_{n-1}^{(\al,\be)}(x).
\label{eq:01.05}
\end{multline}
\begin{multline}
\left(
(2n+\al+\be)(1-x^2)\frac d{dx}
-(n+\al+\be)\bigl((2n+\al+\be)x+\alpha-\beta\bigr)
\right)
P_{n-1}^{(\al,\be)}(x)\\
=-2n(n+\al+\be)\,P_n^{(\al,\be)}(x).
\label{eq:01.09}
\end{multline}
The lowering formula \eqref{eq:01.05} was earlier given in
\cite[10.8(15)]{01}.

In order to obtain \eqref{eq:01.05} from \eqref{eq:01.01},
first rewrite \eqref{eq:01.01} as
\begin{multline*}
\left(
(1-x^2)\frac{d^2}{dx^2}-2(x+\al)\frac d{dx}
\right)
\left(
(1+x)^{-n} P_n^{(\al,\be)}(x)
\right)\\
=-2(n+\al)(n+\be)\,(1+x)^{-n-1} P_{n-1}^{(\al,\be)}(x),
\end{multline*}
and next as
\begin{multline}
\Bigl(
(1+x)(1-x^2)\frac{d^2}{dx^2}+2(1+x)\bigl((n-1)x-n-\al\bigr)\frac d{dx}\\
+n\bigl(-(n-1)x+n+2\al+1\bigr)
\Bigr)P_n^{(\al,\be)}(x)
=-2(n+\al)(n+\be)\,P_{n-1}^{(\al,\be)}(x).
\label{eq:01.03}
\end{multline}
Subtract $(1+x)$ times
the second order differential equation \eqref{eq:01.04}
for Jacobi polynomials from identity \eqref{eq:01.03}
in order to remove its term with a
second order derivative. Then we obtain \eqref{eq:01.05}.

The derivation of \eqref{eq:01.09} from \eqref{eq:01.06} is similar,
with the two intermediate formulas
\begin{multline*}
\left(
(1-x^2)\frac{d^2}{dx^2}-2(x+\al)\frac d{dx}
\right)
\left(
(1+x)^{n+\al+\be} P_{n-1}^{(\al,\be)}(x)
\right)\\
=2n(n+\al+\be)\,(1+x)^{n+\al+\be-1} P_n^{(\al,\be)}(x),
\end{multline*}
\begin{multline}
\Bigl(
(1+x)(1-x^2)\frac{d^2}{dx^2}+2(1+x)\bigl(-(n+\al+\be+1)x+n+\be\bigr)
\frac d{dx}\\
+(n+\al+\be)\bigl(-(n+\al+\be+1)x+n-\al+\be-1\bigr)
\Bigr)
 P_{n-1}^{(\al,\be)}(x)\\
=-2n(n+\al+\be)\,P_n^{(\al,\be)}(x).
\label{eq:01.08}
\end{multline}

The third form of the lowering and raising formulas is:
\begin{multline}
\left(\left(\frac\pa{\pa z}+\frac\pa{\pa w}\right)^2+
\frac{4\be+2}{z+w}\left(\frac\pa{\pa z}+\frac\pa{\pa w}\right)\right)
\left((zw)^n P_n^{(\al,\be)}\left(\thalf(z/w+w/z)\right)\right)\\
=4(n+\al)(n+\be) (zw)^{n-1} P_{n-1}^{(\al,\be)}\left(\thalf(z/w+w/z)\right),
\label{eq:01.12}
\end{multline}
\begin{multline}
\Biggl(\left(z^2\frac\pa{\pa z}+w^2\frac\pa{\pa w}\right)^2+
\left((\al+\be+1)(z+w)-\frac{(2\be+1)zw}{z+w}\right)
\left(z^2\frac\pa{\pa z}+w^2\frac\pa{\pa w}\right)
\\
+(\al+\be+1)\bigl((\al+\be+2)(z^2+w^2)+2(\al-\be)zw\bigr)\Biggr)
\left((zw)^{n-1} P_{n-1}^{(\al,\be)}\left(\thalf(z/w+w/z)\right)\right)
\\
=4n(n+\al+\be)\,(zw)^n P_n^{(\al,\be)}\left(\thalf(z/w+w/z)\right).
\label{eq:01.16}
\end{multline}
The lowering formula \eqref{eq:01.12} can be obtained by rewriting
\eqref{eq:01.03} (use the symmetry \eqref{eq:01.11}).
Similarly, the raising formula \eqref{eq:01.16} is obtained from
\eqref{eq:01.08}.

The cases $\be=\pm\half$ of \eqref{eq:01.12} correspond to iterated cases of
\eqref{eq:02.10} in view of the quadratic transformations
\begin{align}
C_{2n}^{(\la)}(x)=&\frac{(\la)_n}{(\thalf)_n}\,
P_n^{(\la-\half,-\half)}(2x^2-1),
\label{eq.01.13}\\
C_{2n+1}^{(\la)}(x)=&\frac{(\la)_{n+1}}{(\thalf)_{n+1}}\,
x\, P_n^{(\la-\half,\half)}(2x^2-1).
\label{eq.01.14}
\end{align}
\section{Jacobi polynomials for root system $BC_2$}
\label{sec:03}
Jacobi polynomials for root system $BC_2$ are a very special case
(in fact one of the motivating examples) of the Jacobi polynomials
associated with root systems of Heckman and Opdam \cite[Theorem 8.3]{03},
\cite{05}, \cite{04}. (Of course, the $BC_1$ case is given by
the classical Jacobi polynomials of \S\ref{sec:01}.)
They were introduced by the author in
\cite{06}, and further elaborated in
\cite{07} and in (my main reference) \cite{08}.

The \emph{$BC_2$ Jacobi polynomials} $R_{n,k}^{\al,\be,\ga}(\xi,\eta)$
($n\ge k\ge 0$)
are obtained by orthogonalizing the sequence
$1,\xi,\eta,\xi^2,\xi\eta,\eta^2,\xi^3,\ldots,\xi^n,\xi^{n-1}\eta,\ldots,
\xi^{n-k}\eta^k,\ldots$ with respect to the weight function
$\eta^\al(1-\xi+\eta)^\be(\xi^2-4\eta)^\ga\,d\xi\,d\eta$ on the region
in the $(\xi,\eta)$ plane
bounded by the straight lines $\eta=0$ and $1-\xi+\eta=0$ and by
the parabola $\xi^2-4\eta=0$ (so the region has vertices $(0,0)$, $(1.0)$ and
$(2,1)$). Furthermore, the polynomials are normalized such that
$R_{n,k}^{(\al,\be,\ga)}(0,0)=1$.
In fact, it can be shown that $R_{n,k}^{\al,\be,\ga}(\xi,\eta)$ is
only a linear combination of the monomials $\xi^{m-l}\eta^l$
for $m\le n$ and $m+l\le n+k$.

Important special cases of the $BC_2$ Jacobi polynomials occur
for $\ga=\pm\thalf$, where they can be expressed in terms of classical
Jacobi polynomials $R_n^{(\al,\be)}(x)$ (the normalized form, see
\eqref{eq:01.13}):
\begin{multline}
R_{n,k}^{\al,\be,-\half}(x+y,xy)\\
=
\thalf\left(R_n^{(\al,\be)}(1-2x)\,R_k^{(\al,\be)}(1-2y)
+R_k^{(\al,\be)}(1-2x)\,R_n^{(\al,\be)}(1-2y)
\right),
\label{eq:03.03}
\end{multline}
\begin{multline}
R_{n,k}^{\al,\be,\half}(x+y,xy)=
\frac{-(\al+1)}{(n-k+1)(n+k+\al+\be+2)\,(x-y)}\\
\times\left(R_{n+1}^{(\al,\be)}(1-2x)\,R_k^{(\al,\be)}(1-2y)-
R_k^{(\al,\be)}(1-2x)\,R_{n+1}^{(\al,\be)}(1-2y)\right).
\label{eq:03.04}
\end{multline}

Many explicit formulas for $BC_2$ Jacobi polynomials with general values
for the parameters $\al,\be,\ga$ were found in \cite{06}, \cite{07},
\cite{08} by first deriving the desired formula for $\ga=\pm\thalf$,
next guessing the formula for general $\ga$ by interpolation between the
two known cases ($\ga=\pm\thalf$), and finally proving the conjectured
formula in some way. This method was for instance successful in the derivation
of explicit second order differential operators raising or lowering
some parameters
(so-called \emph{shift operators}, which were important
motivating examples for Opdam \cite{05}):
\begin{align*}
D_-^\ga&\colon R_{n,k}^{\al,\be,\ga}\to R_{n-1,k-1}^{\al+1,\be+1,\ga},&
D_+^{\al,\be,\ga}&\colon
R_{n-1,k-1}^{\al+1,\be+1,\ga}\to R_{n,k}^{\al,\be,\ga},\\
E_-^{\al,\be}&\colon R_{n,k}^{\al,\be,\ga}\to R_{n-1,k}^{\al,\be,\ga+1},&
E_+^{\al,\be,\ga}&\colon R_{n-1,k}^{\al,\be,\ga+1}\to R_{n,k}^{\al,\be,\ga}.
\end{align*}
For instance,
\begin{equation*}
D_-^\ga=\tfrac14\left(\pa_{\xi\xi}+\xi\,\pa_{\xi\eta}+\eta\,\pa_{\eta\eta}+
(\ga+\tfrac32)\,\pa_\eta\right).
\end{equation*}
(Here and in the following I use notation $\pa_x$ for the partial derivative
with respect to $x$, and similarly for other variables.)

For \emph{Jack polynomials in two variables} I will use the notation
\begin{equation}
Z_{m,n}^{\ga-\half}(x+y,xy):=P_{m,n}^{1/\ga}(x,y),
\label{eq:03.06}
\end{equation}
where the standard notation for Jack polynomials is used on the
\RHS. These polynomials can be expressed in terms of ultraspherical or
Jacobi polynomials by
\begin{equation}
\begin{split}
Z_{m,n}^{\ga-\half}(x+y,xy)&=
\frac{(m-n)!}{(\ga)_{m-n}}\,(xy)^{\half(m+n)}\,
C_{m-n}^\ga\left(\frac{x+y}{2(xy)^\half}\right)\\
&=\frac{(2\ga)_{m-n}}{(\ga)_{m-n}}\,(xy)^{\half(m+n)}\,
R_{m-n}^{(\ga-\half,\ga-\half)}\left(\frac{x+y}{2(xy)^\half}\right).
\end{split}
\label{eq:03.07}
\end{equation}

The $BC_2$ Jacobi polynomials can be explicitly expanded in terms of
Jack polynomials in two variables (see \cite[Corollary 6.6]{08}):
\begin{equation}
R_{n,k}^{\al,\be,\ga}=
\sum_{l=0}^k \sum_{m=l}^n c_{n,k;m.l}^{\al,\be,\ga}\,Z_{m,l}^\ga,
\label{eq:03.08}
\end{equation}
where
\begin{multline}
c_{n,k;m.l}^{\al,\be,\ga}=
\frac{(-k)_l(-n-\ga-\thalf)_l}{(-n)_l(\al+1)_l}\\
\times\frac{(-n)_m(n+\al+\be+\ga+\tfrac32)_m}
{(\al+\ga+\tfrac32)_m(\ga+\tfrac32)_m}\,
\frac{(k+\al+\be+1)_l(\ga+\tfrac32)_{m-l}}{l!\,(m-l)!}\\
\times\hyp 43{-m+l,-n+k,-n-k-\al-\be-1,\ga+\thalf}
{-n+l,-n-m-\al-\be-\ga-\thalf,2\ga+1}1.
\label{eq:03.02}
\end{multline}
More generally, $BC_n$ Jacobi polynomials can be expanded in terms of
Jack polynomials in $n$ variables with the expansion coefficients given
combinatorially. For this formula, due to Macdonald, see
\cite[(5.12), (5.13)]{09}.
For $k=0$ formulas \eqref{eq:03.02}, \eqref{eq:03.08} simplify to:
\begin{equation}
R_{n,0}^{\al,\be,\ga}=
\sum_{m=0}^n\frac{(-n)_m (n+\al+\be+2\ga+2)_m(\ga+\thalf)_m}
{(\al+\ga+\tfrac32)_m (2\ga+1)_m m!}\,Z_{m,0}^\ga.
\label{eq:03.05}
\end{equation}
\section{Lowering and raising $BC_2$ Jacobi polynomials in a special case}
\label{sec:04}
Let us now look for lowering and raising operators
\[
M_n^{\al,\be,\ga}\colon R_{n,0}^{\al,\be,\ga}\to R_{n-1,0}^{\al,\be,\ga},
\qquad
K_n^{\al,\be,\ga}\colon R_{n,0}^{\al,\be,\ga}\to R_{n+1,0}^{\al,\be,\ga}.
\]
In this paper we will restrict to the case that
$\al=\be$. Let us first try to find such operators acting on
$R_n^{(\al,\al)}(x)\pm R_n^{\al,\al}(y)$
(slight variants of
the case $\al=\be$, $k=0$ of \eqref{eq:03.03} and \eqref{eq:03.04}).
From \eqref{eq:02.06} we obtain
\begin{equation}
\left((1-x^2)\pa_x+nx\right) R_n^{(\al,\al)}(x)=n R_{n-1}^{(\al,\al)}(x).
\label{e	q:04.01}
\end{equation}
Then
\begin{multline*}
\left((1-x^2)\pa_x+nx+(1-y^2)\pa_y+ny\right)
\left(R_n^{(\al,\al)}(x)+R_n^{(\al,\al)}(y)\right)\\
=
n\left(R_{n-1}^{(\al,\al)}(x)+R_{n-1}^{(\al,\al)}(y)\right)+
n\left(x R_n^{(\al,\al)}(y)+y R_n^{(\al,\al)}(x)\right).
\end{multline*}
Here we cannot express all occurrences of $R_m^{(\al,\be)}(x)$ and
$R_m^{(\al,\be)}(x)$ in terms of $R_m^{(\al,\be)}(x)+R_m^{(\al,\be)}(x)$.
The following trick will help us.

Rewrite \eqref{e	q:04.01} as
\begin{equation}
\left((n+2\al+1)(1-x^2)\pa_x+n(n+2\al+1)x\right) R_n^{(\al,\al)}(x)=
n(n+2\al+1) R_{n-1}^{(\al,\al)}(x)
\label{eq:04.02}
\end{equation}
Then recognize $n(n+2\al+1)$ as the eigenvalue in the second order
differential equation for $R_n^{(\al,\al)}(x)$ (see \eqref{eq:01.04}):
\begin{equation}
\left((1-x^2)\pa_{xx}-2(\al+1)x\pa_x\right) R_n^{(\al,\al)}(x)=
-n(n+2\al+1) R_n^{(\al,\al)}(x).
\label{eq:04.03}
\end{equation}
From \eqref{eq:04.02} and \eqref{eq:04.03} we obtain
\begin{multline*}
\left((1-x^2)x\pa_{xx}-
((2\al+2)x^2+(n+2\al+1)(1-x^2))\pa_x\right)R_n^{(\al,\al)}(x)\\
=
-n(n+2\al+1)R_{n-1}^{(\al,\al)}(x).
\end{multline*}
This can be rewritten as
\begin{multline}
\left(x(1-x)(1-2x)\pa_{xx}+(\al+1+2(n-1)x-2(n-1)x^2)\pa_x\right)
R_n^{(\al,\be)}(1-2x)\\
=-n(n+2\al+1) R_{n-1}^{(\al,\al)}(1-2x).
\label{eq:04.05}
\end{multline}
If we add or subtract \eqref{eq:04.05}
and the same identity with $x$ replaced by $y$
then we obtain a lowering operator
acting on $R_n^{(\al,\al)}(1-2x)\pm R_n^{(\al,\al)}(1-2y)$.
There is still a lot of freedom here, since we can add
terms which end on $\pa_{xy}$. Thus, there are many ways to write down
lowering operators acting on $R_{n,0}^{\al,\be,\pm\half}(x+y,xy)$ and it will
be hard to decide in this way on a possible interpolation with respect to
the parameter $\ga$ of the lowering operators for $\ga=\pm\half$.

We can do better by the following approach. Put
\begin{equation}
\FSD_-:=\pa_x+\pa_y,\quad
\FSD_+^\ga:=x^2\pa_x+y^2\pa_y+(\ga+\thalf)(x+y),\quad
\FSD_0:=x\pa_x+y\pa_y.
\label{eq:04.12}
\end{equation}
From \eqref{eq:02.10}, \eqref{eq:02.11} and the homogeneity
of $Z_{m,0}^\ga(x+y,xy)$ in $x,y$ we obtain:
\begin{align}
\FSD_- Z_{m,0}^\ga(x+y,xy)&=
\frac{m(2\ga+m)}{\ga+m-\thalf}\,Z_{m-1,0}^\ga(x+y,xy),
\label{eq:04.06}\\[\smallskipamount]
\FSD_+^\ga Z_{m,0}^\ga(x+y,xy)&=(\ga+m+\thalf) Z_{m+1,0}^\ga(x+y,xy),
\label{eq:04.07}\\[\smallskipamount]
\FSD_0 Z_{m,0}^\ga(x+y,xy)&=
m Z_{m,0}^\ga(x+y,xy).
\label{eq:04.08}
\end{align}
Let us try to use the operators \eqref{eq:04.12} as building blocks
for a lowering operator acting on $R_{n,0}^{\al,\be,\ga}(x+y,xy)$
such that it reduces for $\ga=\pm\thalf$
to an operator we already know. I will work this out here only for the case
$\al=\be$. The following conjectured
lowering formula is obtained:
\begin{multline}
\Bigl(\FSD_-\FSD_0-3(\FSD_0)^2+2\FSD_+^\ga\FSD_0+(\al+\ga+\thalf)\FSD_-
+2(n-2\ga-\thalf)\FSD_0\\
-2n(\FSD_+^\ga-\ga-\thalf)\Bigr)
R_{n,0}^{\al,\al,\ga}(x+y,xy)=-n(n+2\al+1) R_{n-1,0}^{\al,\al,\ga}(x+y,xy).
\label{eq:04.10}
\end{multline}
Formula \eqref{eq:04.10} can indeed be verified
by using \eqref{eq:03.05},
\eqref{eq:04.06}, \eqref{eq:04.07}, \eqref{eq:04.08}.

Similarly as for \eqref{eq:04.10}, one can conjecture
and next prove the following raising formula:
\begin{multline}
\Bigl(\FSD_-\FSD_0-3(\FSD_0)^2+2\FSD_+^\ga\FSD_0+(\al+\ga+\thalf)\FSD_-
-2(n+2\al+4\ga+\tfrac52)\FSD_0\\
+2(n+2\al+2\ga+2)(\FSD_+^\ga-\ga-\thalf)\Bigr)
R_{n,0}^{\al,\al,\ga}(x+y,xy)\\
=-(n+2\ga+1)(n+2\al+2\ga+2) R_{n+1,0}^{\al,\al,\ga}(x+y,xy).
\label{eq:04.11}
\end{multline}
The computations to check \eqref{eq:04.10} and \eqref{eq:04.11}
are feasible on paper, but I have also checked the results
in {\em Mathematica}.
\bibliographystyle{amsalpha}

\end{document}